\documentclass[oneside,notitlepage,12pt]{article}

\usepackage{amssymb}
\usepackage[leqno]{amsmath}
\usepackage{amsfonts}
\usepackage{amsopn}
\usepackage{amstext}
\usepackage{amsthm}
\usepackage{tikz}
\usetikzlibrary{cd}

\usepackage[colorlinks, backref]{hyperref}
\usepackage{calrsfs}
\usepackage{fourier-orns}
\usepackage{clock} %\ClockFrametrue\ClockStyle2
\renewcommand{\Bbb}{\mathbb}

\newenvironment{pf}{\begin{proof}}{\end{proof}}

\newcommand{\Nat}{{\Bbb{N}}}

\newcommand{\lam}{{\lambda}}

\newcommand{\eps}{\varepsilon}
\renewcommand{\phi}{\varphi}
\renewcommand{\rho}{\varrho}

\newcommand{\rest}{\restriction}

\newcommand{\ntr}{{n\in\omega}}
\newcommand{\Ntr}{n\in{\Bbb{N}}}
\newcommand{\loe}{\leq}
\newcommand{\goe}{\geq}

\newcommand{\subs}{\subseteq}
\newcommand{\sups}{\supseteq}

\newcommand{\dist}{\operatorname{dist}}

\newcommand{\id}[1]{{\operatorname{i\!d}_{#1}}} % identity morphism

\newcommand{\oraz}{\qquad\text{and}\qquad}

\newtheorem{tw}{Theorem}[section]
\newtheorem{wn}[tw]{Corollary}

\newtheorem{lm}[tw]{Lemma}
\newtheorem{prop}[tw]{Proposition}
\newtheorem{claim}[tw]{Claim}
\theoremstyle{definition}

\newtheorem{ex}[tw]{Example}
\newtheorem{pyt}[tw]{Question}

\theoremstyle{remark}

\newcommand{\setof}[2]{\{#1\colon #2\}}
\newcommand{\Bigsetof}[2]{\Bigl\{#1\colon #2\Bigr\}}

\newcommand{\sett}[2]{\{#1\}_{#2}}
\newcommand{\sn}[1]{\{#1\}} % singleton
\newcommand{\map}[3]{#1\colon #2 \to #3} % A function
\newcommand{\img}[2]{#1[#2]} % image of a set
\newcommand{\ciag}[1]{{\sett{{#1}_n}{\ntr}}}

\newcommand{\anorm}{\|\cdot\|}
\newcommand{\norm}[1]{\|#1\|}

\newcommand{\bS}{{\mathbb{S}}}

\newcommand{\cmp}{\circ} % composition!

\newcommand{\C}{{\ensuremath\mathcal C}} % The "continuous functions" functor.

\newcommand{\separator}{\begin{center}***\end{center}}

\newcommand{\N}{\mathbb N}

\newcommand{\G}{\mathbb G}

\newcommand{\intoc}[2]{(#1, #2]}

\newcommand{\Gurarii}{Gurari\u\i}
\newcommand{\uop}{\mathbf\Omega} % The universal operator.
\newcommand{\upro}[1]{\mathbf P_{#1}}

\title{A note on universal operators between separable Banach spaces}
\author{
{\sc Joanna Garbuli\'nska-W\c egrzyn}\\
{\small Institute of Mathematics,}
{\small Jan Kochanowski University, Poland}
\and
{\sc Wies{\l}aw Kubi\'s}\\
{\small Institute of Mathematics, Czech Academy of Sciences, Czech Republic}\\
{\small Department of Mathematics, Cardinal Stefan Wyszy\'nski University in Warsaw}
}

\begin{document}
\maketitle

\begin{abstract}
We compare two types of universal operators constructed relatively recently by Cabello S\'anchez, and the authors.
The first operator $\uop$ acts on the \Gurarii\ space, while the second one $\upro \bS$ has values in a fixed separable Banach space $\bS$.
We show that if $\bS$ is the \Gurarii\ space, then both operators are isometric.
We also prove that, for a fixed space $\bS$, the operator $\upro{\bS}$ is isometrically unique.
Finally, we show that 
$\uop$ is {generic} in the sense of a natural infinite game.

\ 


\end{abstract}

\tableofcontents

\section{Universal operators}

The purpose of this note is to discuss two constructions of universal operators between separable Banach spaces.
We are interested in isometric universality.
Namely, an operator $U$ is \emph{universal}
if its restrictions to closed subspaces are, up to linear isometries, \emph{all} linear operators of norm not exceeding $\norm{U}$.
To be more precise, a bounded linear operator $\map{U}{V}{W}$ acting between separable Banach spaces is \emph{universal} if for every linear operator $\map{T}{X}{Y}$ with $X$, $Y$ separable and $\norm{T} \loe \norm{U}$, there exist linear isometric embeddings
$\map i X V$, $\map j Y W$
such that the diagram
$$\begin{tikzcd}
	V \ar[rr, "U"] & & W \\
	X \ar[rr, "T"']\ar[u, hook, "i"] & & Y \ar[u, hook, "j"']
\end{tikzcd}$$
is commutative, that is, $U \cmp i = j \cmp T$.
Such an operator has been relatively recently constructed by the authors~\cite{GK}.
Another recent work~\cite{CsGwK}, due to Cabello S\'anchez and the present authors, contains in particular a construction of a linear operator that is universal in a different sense.
Namely, let us say that a bounded linear operator $\map{U}{V}{W}$ is \emph{left-universal} (for operators into $W$) if for every linear operator $\map T X W$ with $X$ separable and $\norm{T} \loe \norm{U}$ there exists a linear isometric embedding $\map{i}{X}{V}$ for which the diagram
$$\begin{tikzcd}
	V \ar[rr, "U"] & & W \\
	X \ar[rru, "T"'] \ar[u, hook, "i"] & &
\end{tikzcd}$$
is commutative, that is, $U \cmp i = T$.
Clearly, if $W$ is isometrically universal in the class of all separable Banach spaces then a left-universal operator with values into $W$ is universal.
The left-universal operator $U$ constructed in \cite{CsGwK} had been later essentially used (with a suitable space $W$) for finding an isometrically universal graded Fr\'echet space~\cite{BarKakKub}.
There exist other concepts of universality in operator theory, see the introduction of \cite{GK} for more details and references.

Let us note the following simple facts related to universal operators.

\begin{prop}
	Let $\map U V W$ be a bounded linear operator acting between separable Banach spaces.
	\begin{enumerate}
		\item[$(1)$] If $U$ is universal then both $V$ and $W$ are isometrically universal among the class of separable Banach spaces.
		\item[$(2)$] Assume $U$ is left-universal.
		Then $\ker U$ is isometrically universal among the class of separable Banach spaces.
		Furthermore, $U$ is right-invertible, that is, there exists an isometric embedding $\map e W V$ such that $U \cmp e = \id{W}$.
		\item[$(3)$] Assume $U$ is (left-)universal. Then $\lam U$ is (left-)universal for every $\lam>0$.
	\end{enumerate}
\end{prop}

\begin{pf}
	(1) Fix a separable Banach space $X$. Taking the zero operator $\map T X 0$, we see that $V$ contains an isometric copy of $X$. Taking the identity $\id X$, we see that $W$ contains an isometric copy of $X$.
	
	(2) The same argument as above, using the zero operators, shows that $\ker U$ is isometrically universal.
	Taking the identity $\id{W}$, we obtain the required isometric embedding $\map e W V$.
	
	(3) Assume $U$ is universal, fix $\lam>0$ and fix $\map T X Y$ with $\norm T \loe \lam \norm U$.
	Then $\norm{\lam^{-1}T} \loe \norm U$, therefore there are isometric embeddings $\map i X V$, $\map j Y W$ such that $U \cmp i = j \cmp (\lam^{-1} T)$.
	Finally, $(\lam U) \cmp i = j \cmp T$.
	If $U$ is left-universal, the argument is the same, the only difference is that $j = \id W$.
\end{pf}

By (3) above, we may restrict attention to non-expansive operators.
It turns out that there is an easy way of constructing left-universal operators, once we have in hand an isometrically universal space.
The argument below was pointed out to us by Przemys{\l}aw Wojtaszczyk.

\begin{ex}\label{ExUnivOptLeft}
	Let $V$ be an isometrically universal Banach space and let $W$ be an arbitrary Banach space.
	Consider $V \oplus W$ with the maximum norm and let $$\map{\pi}{V \oplus W}{W}$$ be the canonical projection.
	Given a non-expansive operator $\map T X W$ with $X$ separable, choose an isometric embedding $\map e X V$ and define $\map{j}{X}{V \oplus W}$ by $j(x) = (e(x), T(x))$.
	Then $j$ is an isometric embedding and $\pi \cmp j = T$, showing that $\pi$ is left-universal.
	Of course, if additionally $W$ is isometrically universal, then $\pi$ is a universal operator.	
\end{ex}

Perhaps the most well known universal Banach space is $\C([0,1])$, the space of all continuous (real or complex) valued functions on the unit interval, endowed with the maximum norm. In view of the example above, there exists a universal operator from $\C([0,1]) \oplus \C([0,1])$ onto $\C([0,1])$.
This leads to (at least potentially) many other universal operators, namely:

\begin{prop}\label{PROPsdkvbdag}
	Let $V, W$ be isometrically universal separable Banach spaces.
	Then there exists a universal operator from $V$ into $W$.
\end{prop}

\begin{pf}

	Fix a universal operator $\map{\pi}{E}{F}$ (for instance, $E = \C([0,1]) \oplus \C([0,1])$ and $F = \C([0,1])$) and fix a linear isometric embedding $\map e E V$.
	Using the amalgamation property for linear operators, we find a separable Banach space $V'$, a linear isometric embedding $\map{e'}{F}{V'}$, and a non-expansive linear operator $\map{\Omega}{V}{V'}$ for which the diagram
	$$\begin{tikzcd}
	V \ar[rr, "{\Omega}"] & & V' \\
	E \ar[rr, "\pi"'] \ar[u, hook, "e"] & & F \ar[u, hook, "e'"']
	\end{tikzcd}$$
	is commutative.
	As $W$ is isometrically universal, we may additionally assume that $V' = W$, replacing $\Omega$ by $i \cmp \Omega$ and $e'$ by $i \cmp e'$, where $i$ is a fixed isometric embedding of $V'$ into $W$.
	It is evident that now $\Omega$ is a universal operator, because of the universality of $\pi$.
\end{pf}

As a consequence, there exists a universal operator on $\C([0,1])$.
We do not know whether there exists a left-universal operator on $\C([0,1])$. The situation changes when replacing $[0,1]$ with the Cantor set $2^\Nat$. The space $\C(2^\Nat)$ is linearly isomorphic (but not isometric) to $\C([0,1])$ and it is isometrically universal, too.
Furthermore, $\C(2^\Nat) \oplus \C(2^\Nat)$ with the maximum norm is linearly isometric to $\C(2^\Nat)$, because the disjoint sum of two copies of the Cantor set is homeomorphic to the Cantor set. Thus, Example~\ref{ExUnivOptLeft} provides a left-universal operator on $\C(2^\Nat)$.

Another, not so well known, universal Banach space is the \emph{\Gurarii\ space}.
This is the unique, up to a linear isometry, separable Banach space $\G$ satisfying the following condition:
\begin{enumerate}
	\item[(G)] For every $\eps>0$, for every finite-dimensional spaces $X_0 \subs X$, for every linear isometric embedding $\map{f_0}{X_0}{\G}$ there exists a linear $\eps$-isometric embedding $\map{f}{X}{\G}$ such that $f \rest X_0 = f_0$.
\end{enumerate}
By an \emph{$\eps$-isometric embedding} (briefly: \emph{$\eps$-embedding}) we mean a linear operator $f$ satisfying
$$(1-\eps)\norm x \loe \norm{f(x)} \loe (1+\eps)\norm x$$
for every $x$ in the domain of $f$.
The space $\G$ was constructed by \Gurarii~\cite{gurarii}; its uniqueness was proved by Lusky~\cite{lusky}.

The universal operator constructed in~\cite{GK} has a special property that actually makes it unique, up to linear isometries.
Below we quote the precise result.
 
\begin{tw}[{\cite{GK}}]\label{GK}
	There exists a non-expansive linear operator $\map{\uop}{\G}{\G}$ with the following property:
	\begin{enumerate}
		\item[{\rm(G)}] Given $\eps>0$, given a non-expansive operator $\map {T}{X}{Y}$ between finite-dimensional spaces, given $X_0 \subs X$, $Y_0 \subs Y$ and isometric embeddings $\map i {X_0}U$, $\map j {Y_0}V$ such that $\uop \cmp i = j \cmp (T\rest X_0)$, there exist $\eps$-embeddings $\map {i'} X U$, $\map {j'} Y V$ satisfying
		$$\norm{i'\rest X_0 - i} \loe \eps, \quad \norm{j'\rest Y_0 - j} \loe \eps, \oraz \norm{\uop \cmp i'  - j' \cmp T} \loe \eps.$$
	\end{enumerate}
	Furthermore, $\uop$ is a universal operator and property $(G)$ specifies it uniquely, up to a linear isometry.
\end{tw}
According to \cite{GK}, we shall call condition (G) the \emph{\Gurarii\ property}.
%We shall later see that the \Gurarii\ property implies a stronger one, without inequalities.
What makes this operator of particular interest is perhaps its \emph{almost homogeneity}:

\begin{tw}[{\cite{GK}}]
	Given finite-dimensional subspaces $X_0,X_1,Y_0,Y_1$ of $\G$, given linear isometries $\map i{X_0}{X_1}$, $\map j{Y_0}{Y_1}$ such that $\uop \cmp i = j \cmp \uop$, for every $\eps>0$ there exist bijective linear isometries $\map I \G \G$, $\map J \G \G$ satisfying $\uop \cmp I = J \cmp \uop$ and $\norm{I \rest X_0 - i} < \eps$, $\norm{J \rest Y_0 - j} < \eps$.
\end{tw}

We now describe the left-universal operators constructed in~\cite{CsGwK}.
Fix a separable Banach space $\bS$.

\begin{tw}[{\cite[Section 6]{CsGwK}}]\label{Thmsibvisbe}
	There exists a non-expansive linear operator $\map{\upro{\bS}}{V_\bS}{\bS}$ with $V_\bS$ a separable Banach space, satisfying the following condition:
	\begin{enumerate}
		\item[$(\ddagger)$] For every finite-dimensional spaces $X_0 \subs X$, for every non-expansive linear operator $\map T X \bS$, for every linear isometric embedding $\map e{X_0}{V_\bS}$ such that $\upro \bS \cmp e = T \rest X_0$, for every $\eps>0$ there exists an $\eps$-embedding $\map f X {V_\bS}$ satisfying
		$$\norm{f \rest X_0 - e} \loe \eps \qquad \text{and} \qquad \norm{\upro \bS \cmp f - T} \loe \eps.$$
	\end{enumerate}
	Furthermore, $\upro \bS$ is left-universal for operators into $\bS$.
\end{tw}

We shall say that an operator $P$ has the \emph{left-\Gurarii} property if it satisfies ($\ddagger$) in place of $\upro \bS$. Of course, unlike the \Gurarii\ property, the left-\Gurarii\ property involves a parameter $\bS$, namely, the common range of the operators.

Actually, the projection $\upro{\bS}$ was constructed in \cite{CsGwK} in case where $\bS$ had some additional property, needed only for determining the domain of $\upro{\bS}$.
Moreover, \cite{CsGwK} deals with $p$-Banach spaces, where $p \in \intoc 01$, however $p=1$ gives exactly the result stated above.
Operators $\upro{\bS}$ have the following property which can be called \emph{almost left-homogeneity}.

\begin{tw}\label{THMalleftHomnity}
Given finite-dimensional subspaces $X_0, X_1$ of $V_\bS$, a linear isometry $\map{h}{X_0}{X_1}$ such that $\upro{\bS} \cmp h = \upro{\bS} \rest X_0$, for every $\eps>0$ there exists a bijective linear isometry $\map H {V_\bS}{\bS}$ satisfying $\upro{\bS} \cmp H = \upro{\bS}$ and $\norm{H \rest X_0 - h} < \eps$.
\end{tw}

In this note we present a proof that condition $(\ddagger)$ determines $\upro{\bS}$ uniquely, up to linear isometries. The arguments will also provide a proof of Theorem~\ref{THMalleftHomnity}. Furthermore, we show that
$\uop = \upro \G$ and that $\uop$ is a generic operator in the space of all non-expansive operators on the \Gurarii\ space into itself, in the sense of a natural variant of the Banach-Mazur game.

\section{Properties of $\uop$ and $\upro{\bS}$}

Let us recall the following easy fact concerning finite-dimensional normed spaces (cf.~\cite[Thm. 2.7]{kubis-gar} or \cite[Claim 2.3]{BarKakKub}).
It actually says that the strong operator topology is equivalent to the norm topology in the space of linear operators with a fixed finite-dimensional domain.

\begin{lm}\label{Lmsdsdbgiwbi}
	Let $A$ be a vector basis of a finite-dimensional normed space $E$.
	For every $\eps>0$ there exists $\delta>0$ such that for every Banach space $X$, for every linear operator $\map{f}{E}{X}$ the following implication holds:
	$$\max_{a \in A} \norm{f(a)}\loe \delta \implies \norm{f}\loe \eps.$$
\end{lm}

\begin{pf}
	Fix $M>0$ satisfying the following condition:
	\begin{enumerate}
	\item[(*)] $\max_{a\in A}|\lam_a| \loe M$ whenever $x = \sum_{a \in A}\lam_a a$ and $\norm x \loe 1$.
	\end{enumerate}
	Such $M$ clearly exists, because of compactness of the unit ball of $E$.
	Now, given $\eps>0$, let $\delta = \eps/(M\cdot |A|)$.
	Suppose $\max_{a\in A}\norm{f(a)}\loe \delta$.
	Then, given $x = \sum_{a \in A}\lam_a a$ with $\norm x \loe 1$, we have
	$$\norm{f(x)} \loe \sum_{a \in A}|\lam_a| \cdot \norm{f(a)} \loe |A| \cdot M \cdot \delta = \eps.$$
	We conclude that $\norm{f} \loe \eps$.
\end{pf}

The following result, in case $\bS = \G$ can be found in~\cite{BarKakKub}.

\begin{tw}\label{ThmIidfbisa}
	Let $\map P V \bS$ be a linear operator.
	The following conditions are equivalent.
	\begin{enumerate}
		\item[\rm{(a)}] $P$ has the left-\Gurarii\ property {\rm($\ddagger$)}.
		\item[\rm{(b)}] For every finite-dimensional spaces $X_0 \subs X$, for every non-expansive linear operator $\map T X \bS$, for every linear isometric embedding $\map e{X_0}{V}$ such that $P \cmp e = T \rest X_0$, for every $\eps>0$ there exists an $\eps$-embedding $\map f X {V}$ satisfying
		$${f \rest X_0 = e} \qquad \text{and} \qquad {P \cmp f = T}.$$
	\end{enumerate}
\end{tw}

\begin{pf}
	Obviously, (b) is stronger than ($\ddagger$).
	
	Fix $\eps>0$ and fix a vector basis $A$ of $X$ such that $A_0 = X_0 \cap A$ is a basis of $X_0$.
	We may assume that $\norm{a}=1$ for every $a \in A$.
	Fix $\delta>0$ and apply the left-\Gurarii\ property for $\delta$ instead of $\eps$.
	We obtain a $\delta$-embedding $\map{f}{X}{V}$ such that $\norm{f\rest X_0 - e} \loe \delta$ and $\norm{P \cmp f - T} \loe \delta$.
	Define $\map{f'}{X}{V}$ by the conditions $f'(a) = e(a)$ for $a \in A_0$ and $f'(a) = f(a)$ for $a \in A \setminus A_0$.
	Note that $\norm{f'(a)-f(a)} \loe \delta$ for every $a \in A$.
	Thus, if $\delta$ is small enough, then by Lemma~\ref{Lmsdsdbgiwbi}, we can obtain that $f'$ is an $\eps$-embedding.
	Furthermore, $\norm{P \cmp f' - P \cmp f} \loe \eps$ (recall that $\delta$ depends on $\eps$ and the norm of $X$ only), therefore $\norm{P \cmp f' - T} \loe \eps + \delta$.
	
	The arguments above show that for every $\eps>0$ there exists an $\eps$-embedding $\map{f'}{X}{V}$ extending $e$ and satisfying $\norm{P \cmp f' - T} \loe \eps$.
	
	Let us apply this property for $\delta$ instead of $\eps$, where $\delta$ is taken from Lemma~\ref{Lmsdsdbgiwbi}.
	We obtain a $\delta$-embedding $\map{f}{X}{V}$ extending $e$ and satisfying $\norm{P \cmp f - T} \loe \delta$.
	
	Given $a \in A \setminus A_0$, the vector
	$$w_a = P(f(a)) - T(a)$$
	has norm $\loe \delta$.
	Define $\map{f'}{X}{V}$ by the conditions $f'\rest X_0 = e$ and
	$$f'(a) = f(a) - w_a$$
	for $a \in A \setminus A_0$.
	Lemma~\ref{Lmsdsdbgiwbi} implies that $f'$ is an $\eps$-embedding, because $\norm{f'(a) - f(a)} = \norm{w_a} \loe \delta$ for $a \in A \setminus A_0$.
	Finally, given $a \in A \setminus A_0$, we have
	$$P f' (a) = P f (a) - w_a = T(a)$$
	and the same obviously holds for $a \in A_0$.
	Thus $P \cmp f' = T$.
\end{pf}

The proof of the next result is just a suitable adaptation of the arguments above, therefore we skip it. 

\begin{prop}
	Let $\map \Omega V W$ be a linear operator. The following conditions are equivalent.
	\begin{enumerate}
		\item[\rm{(a)}] $\Omega$ has the \Gurarii\ property {\rm(G)}.
		\item[\rm{(b)}] Given $\eps>0$, given a non-expansive operator $\map {T}{X}{Y}$ between finite-dimensional spaces, given $X_0 \subs X$, $Y_0 \subs Y$ and isometric embeddings $\map {i_0} {X_0}V$, $\map {j_0}{Y_0}W$ such that $\Omega \cmp i_0 = j_0 \cmp (T\rest X_0)$, there exist $\eps$-embeddings $\map {i} X V$, $\map {j} Y W$ satisfying
		$$i\rest X_0 = i_0, \quad j\rest Y_0 = j_0, \oraz \Omega \cmp i = j \cmp T.$$
	\end{enumerate}
\end{prop}

The last result of this section is the key step towards identifying $\uop$ with $\upro{\G}$.

\begin{tw}\label{ThmLeftuydsn}
	The operator $\uop$ has the left-\Gurarii\ property (i.e., it satisfies condition $(\ddagger)$ of Theorem~\ref{Thmsibvisbe} with $\bS = \G$).
	In particular, it is left-universal.
\end{tw}

\begin{pf}
	Fix a non-expansive linear operator $\map T X \G$ with $X$ finite-dimensional, and fix an isometric embedding $\map{e}{X_0}{\G}$, where $X_0$ is a linear subspace of $X$ and $T\rest X_0 = \uop \cmp e$.
	Let $Y_0 = Y = \img T X \subs \G$ and consider $T$ as an operator from $X$ to $Y$.
	Applying the \Gurarii\ property with $i = e$ and $j$ the inclusion $Y_0 \subs \G$, we obtain an $\eps$-embedding $\map{e'}{X}{\G}$ which is $\eps$-close to $e$ and satisfies $\norm{\uop \cmp e' - T} \loe \eps$.
	This is precisely condition ($\ddagger$) from Theorem~\ref{Thmsibvisbe}.
\end{pf}

In order to conclude that $\uop = \upro{\G}$, it remains to show that ($\ddagger$) determines the operator uniquely.
This is done in the next section.

\section{Uniqueness of $\upro{\bS}$}

Before proving that the left-\Gurarii\ property determines the operator uniquely, we quote the following crucial lemma from~\cite{GarFDD}.

\begin{lm}\label{LmDFAfefqo}
	Let $\eps>0$ and let $\map{f}{E}{F}$ be an $\eps$-embedding, where $E$, $F$ are Banach spaces. Let $\map{\pi}{E}{\bS}$, $\map{\rho}{F}{\bS}$ be non-expansive linear operators such that $\norm{\rho \cmp f - \pi} \loe \eps$.
	Then there exists a norm on $Z = X \oplus Y$ such that the canonical embeddings $\map{i}{X}{Z}$, $\map{j}{Y}{Z}$ are isometric, $\norm{j \cmp f - i} \loe \eps$ and the operator $\map t Z \bS$ defined by
	$t(x,y) = \pi(x) + \rho(y)$
	is non-expansive.
\end{lm}

Note that the operator $t$ satisfies $t \cmp i = \pi$ and $t \cmp j = \rho$.
Actually, the norm mentioned in the lemma above does not depend on the operators $\pi$, $\rho$.
It is defined by the following formula:
\begin{equation*}
\norm{(x,y)} = \inf \Bigsetof{\norm{x-w}_X+\norm{y-f(w)}_Y+\eps\norm{w}_X}{w \in X},
\tag{$*$}\label{EqJTGsgij}
\end{equation*}
where $\anorm_X$, $\anorm_Y$ denote the norm of $X$ and $Y$, respectively.
An easy exercise shows that (\ref{EqJTGsgij}) is the required norm, proving Lemma~\ref{LmDFAfefqo}.

\begin{tw}\label{ThmUniqunsmnn}
	Let $\bS$ be a separable Banach space and let $\map \pi E \bS$, $\map{\pi'}{E'}{\bS}$ be non-expansive linear operators, both with the left-\Gurarii\ property.
	If $E$, $E'$ are separable Banach spaces, then there exists a linear isometry $\map i E{E'}$ such that $\pi = \pi' \cmp i$.
	In particular, $\pi$ and $\pi'$ are linearly isometric to $\upro{\bS}$.
\end{tw}

\begin{pf}
	It suffices to prove the following
	
	\begin{claim}\label{ClaimRbirfbir}
		Let $E_0 \subs E$ be a finite-dimensional space, $0<\eps<1$, let $\map{i_0}{E_0}{E'}$ be an $\eps$-embedding such that $\pi' \cmp i_0 = \pi \rest E_0$.
		Then for every $v \in E$, $v' \in E'$, for every $\eta>0$ there exists an $\eta$-embedding $\map{i_1}{E_1}{E'}$ with $E_1$ finite-dimensional and the following conditions are satisfied:
		\begin{enumerate}
			\item[$(1)$] $v \in E_1$ and $\dist(v',\img{i_1}{E_1}) < \eta$;
			\item[$(2)$] $\norm{i_0 - i_1 \rest E_0} < \eps + \eta$ and $\pi' \cmp i_1 = \pi$.
		\end{enumerate}		
	\end{claim}
	
	Using Claim~\ref{ClaimRbirfbir} together with the separability of $E$ and $E'$, we can construct a sequence $\map{i_n}{E_n}{E'}$ of linear operators such that $i_n$ is a $2^{-n}$-embedding, $\bigcup_\ntr E_n$ is dense in $E$ and $\bigcup_\ntr \img{i_n}{E_n}$ is dense in $E'$ and
	$$\norm{i_n - i_{n+1} \rest E_n} \loe 2^{-n} + 2^{-n-1} \oraz \pi' \cmp i_{n+1} = \pi$$
	for every $\ntr$.
	It is evident that $\ciag{i}$ converges pointwise to a linear isometry whose completion $i$ is the required bijection from $E$ onto $E'$ satisfying $\pi' \cmp i = \pi$.
	Thus, it remains to prove Claim~\ref{ClaimRbirfbir}.
	
	This will be carried out by making two applications of Lemma~\ref{LmDFAfefqo}.
	
	Fix $0<\delta<1$, more precise estimations for $\delta$ will be given later.
	Let $E_0' \subs E'$ be a finite-dimensional space containing $v'$ and such that $\img{i_0}{E_0} \subs E_0'$.
	Applying Lemma~\ref{LmDFAfefqo}, we obtain linear isometric embeddings $\map{e_1}{E_0}{W_0}$, $\map{f_1}{E_0'}{W_0}$ and a non-expansive operator $\map{t_0}{W_0}{\bS}$ such that $t_0 \cmp e_1 = \pi \rest E_0$, $t_0 \cmp f_1 = \pi' \rest E_0'$, and $\norm{e_1 - f_1 \cmp i_0} \loe \eps$.
	Knowing that $\pi$ has the left-\Gurarii\ property, by Theorem~\ref{ThmIidfbisa} applied to the isometric embedding $e_1$, we obtain a $\delta$-embedding $\map{g_1}{W_0}{E}$ such that $g_1 \cmp e_1$ is identity on $E_0$ and $\pi \cmp g_1 = t_0$.
	
	Now note that $g_1 \cmp f_1$ is a $\delta$-embedding of $E_0'$ into a finite-dimensional subspace $E_1$ of $E$. Without loss of generality, we may assume that $v \in E_1$.
	Applying Lemma~\ref{LmDFAfefqo} again to $g_1 \cmp f_1$, we obtain linear isometric embeddings $\map{e_2}{E_1}{W_1}$, $\map{f_2}{E_0'}{W_1}$ and a non-expansive linear operator $\map{t_1}{W_1}{\bS}$ such that $t_1 \cmp e_2 = \pi \rest E_1$, $t_1 \cmp f_2 = \pi' \rest E_0'$, and $\norm{e_2 \cmp g_1 \cmp f_1 - f_2} \loe \delta$.
	Knowing that $\pi'$ has the left-\Gurarii\ property and using Theorem~\ref{Thmsibvisbe} for the isometric embedding $f_2$, we obtain a $\delta$-embedding $\map{g_2}{W_1}{E'}$ such that $g_2 \cmp f_2$ is identity on $E_0'$ and $\pi' \cmp g_2 = t_1$.
	The configuration is described in the following diagram, where the horizontal arrows are inclusions, the triangle $E_0 E_0' W_0$ is $\eps$-commutative, and the triangle $E_0' E_1 W_1$ is $\delta$-commutative.
	$$\begin{tikzcd}
			E_0 \ar[rr] \ar[dd, "i_0"'] \ar[dr, "e_1"] & & E_1 \ar[rrrr] \ar[dr, "e_2"] & & & & E \\
			& W_0 \ar[ur, "g_1"'] & & W_1 \ar[drrr, "g_2"] \\
			E_0' \ar[rrrrrr] \ar[ur, "f_1"] \ar[urrr, "f_2"] & & & & & & E'		
	\end{tikzcd}$$	
	It remains to check that $i_1 := g_2 \cmp e_2$ is the required $\delta$-embedding.
	
	First, recall that $v \in E_1$, $v' \in E_0'$ and $v' = g_2(f_2(v'))$.
	Thus, using the fact that $\norm{g_2} \loe 1+\delta$, we get
	\begin{align*}
		\norm{i_1 g_1 f_1(v') - v'} &= \norm{g_2 e_2 g_1 f_1(v') - g_2 f_2(v')}\\
		&\loe (1+\delta) \norm{e_2 g_1 f_1(v') - f_2(v')}\\
		& \loe (1+\delta)\delta \norm{v'}.
	\end{align*}
	Now if $(1+\delta)\delta \norm{v'}<\eta$, then we conclude that $\dist(v',\img{i_1}{E_1}) < \eta$, therefore condition (1) is satisfied.
	
	Given $x \in E_1$, note that
	\begin{align*}
		\pi' i_1(x) = \pi' g_2 e_2 (x) = t_1 e_2(x) = \pi(x).
	\end{align*}
	Here we have used the fact that $\pi' \cmp g_2 = t_1$ and $t_ \cmp e_2 = \pi \rest E_1$.
	
	Furthermore, given $x \in E_0$, we have
	\begin{align*}
		\norm{i_1(x) - i_0(x)} &= \norm{g_2 e_2(x) - i_0(x)} = \norm{g_2 e_2 g_1 e_1(x) - g_2 f_2 i_0(x)} \\
		&\loe (1+\delta) \norm{e_2 g_1 e_1(x) - f_2 i_0(x)},
	\end{align*}
	because $\norm{g_2} \loe 1+\delta$.
	On the other hand,
	\begin{align*}
		\norm{e_2 g_1 e_1(x) - f_2 i_0(x)} &\loe \norm{e_2 g_1 e_1(x) - e_2 g_1 f_1 i_0(x)} + \norm{e_2 g_1 f_1 i_0(x) - f_2 i_0(x)} \\
		&= \norm{g_1 e_1(x) - g_1 f_1 i_0(x)} + \norm{e_2 g_1 f_1 i_0(x) - f_2 i_0(x)} \\
		&\loe (1+\delta)\norm{e_1(x) - f_1 i_0(x)} + \delta \norm{i_0(x)} \\
		&\loe (1+\delta) \eps \norm x + \delta (1+\eps) \norm x \loe (\eps+3\delta) \norm x.
	\end{align*}
	Here we have used the following facts: $e_2$ is an isometric embedding, $g_1$ is a $\delta$-embedding, $i_0$ is an $\eps$-embedding, $\norm{e_2 g_1 f_1 - f_2} \loe \delta$, $\norm{e_1 - f_1 i_0} \loe \eps$ and $\eps<1$.
	
	Finally, $\norm{i_1(x) - i_0(x)} \loe (1+\delta) (\eps + 3\delta) \norm x \loe (\eps + 7\delta) \norm x$.
	Summarizing, if $(1+\delta)\delta \norm{v'}<\eta$ and $7\delta < \eta$ then conditions (1), (2) are satisfied.
	This completes the proof.
\end{pf}

Note that if $\bS$ is the trivial space, the proof above reduces to the well known uniqueness of the \Gurarii\ space, shown by this way in~\cite{KS}.
Furthermore, the arguments above can be applied to $\pi = \pi' = \upro{\bS}$ and $i_0 = h$, thus proving Theorem~\ref{THMalleftHomnity}.
Theorems~\ref{ThmLeftuydsn} and~\ref{ThmUniqunsmnn} yield the following result, announced before.

\begin{wn}
	$\uop = \upro{\G}$.
\end{wn}

In particular, $V_\G = \G$.
It has been shown in~\cite{CsGwK} that $V_\bS = \G$ as long as $\bS$ is a (separable) \emph{Lindenstrauss space}, namely, an isometric $L_1$ predual or (equivalently) a locally almost 1-injective space. Instead of going into details, let us just say that Lindenstrauss spaces are those (separable) Banach spaces that are linearly isometric to a 1-complemented subspace of the \Gurarii\ space. The non-trivial direction was proved by Wojtaszczyk~\cite{Wojt}.
Thus, since $\upro{\bS}$ is a projection, if $V_\bS$ is linearly isometric to $\G$ then $\bS$ is necessarily a Lindenstrauss space.

\section{Generic operators}

Inspired by the result of~\cite{kubBMGurarii}, let us consider the following infinite game for two players \emph{Eve} and \emph{Adam}.
Namely, Eve starts by choosing a non-expansive linear operator $\map{T_0}{E_0}{F_0}$, where $E_0$, $F_0$ are finite-dimensional normed spaces.
Adam responds by a non-expansive linear operator $\map{T_1}{E_1}{F_1}$, such that $E_1 \sups E_0$, $F_1 \sups F_0$ are again finite-dimensional and $T_1$ extends $T_0$.
Eve responds by a further non-expansive linear extension $\map{T_2}{E_2}{F_2}$, and so on.
So at each stage of the game we have a linear operator between finite-dimensional normed spaces.
After infinitely many steps we obtain a chain of non-expansive operators $\sett{\map{T_n}{E_n}{F_n}}{\ntr}$. Let $\map{T_\infty}{E_\infty}{F_\infty}$ denote the completion of its union, namely, $E_\infty$ is the completion of $\sett{E_n}{\ntr}$, $F_\infty$ is the completion of $\sett{F_n}{\ntr}$ and $T_\infty \rest E_n = T_n$ for every $\ntr$.
So far, we cannot say who wins the game.

Let us say that a (necessarily non-expansive) linear operator $\map U X Y$ is \emph{generic} if Adam has a strategy making the operator $T_\infty$ isometric to $U$.
Recall that operators $U,V$ are \emph{isometric} if there are bijective linear isometries $i$, $j$ such that $U \cmp j = i \cmp V$.

\begin{tw}\label{THMsztyryJedna}
	The operator $\uop$ is generic.
\end{tw} 

\begin{pf} Let us fix a non-expansive linear operator $U:\G \to \G$ between separable Banach spaces satisfying (G). 
Adam’s strategy can be described as follows.

Fix a countable set $\{v_n: a_n\to b_n\}_{n\in \N}$ linearly dense in $U:\G \to \G$. Let $\map{T_0}{E_0}{F_0}$ be the first move of Eve. Adam finds isometric embeddings $i_0:E_0\to \G$, $j_0:F_0\to \G$ and  finite-dimensional spaces $E_0\subset E_1$, $F_0 \subset F_1$ together with isometric embeddings $i_1:E_1\to \G$ , $j_1:F_1\to \G$ and non-expansive linear operators $\map{T_1}{E_1}{F_1}$ such that $T_1$ extends $T_0$, $a_0\in i_1[E_1]$, $b_0\in j_1[F_1]$.

Suppose now that $n = 2k > 0$ and $T_n:E_n \to F_n$ was the last move of Eve. We assume that
linear isometric embeddings $i_{n-1} : E_{n-1}\to \G$,  $j_{n-1} : F_{n-1}\to \G$ have already been fixed. Using (G) from Theorem \ref{GK} we choose
linear isometric embeddings $i_{n} : E_{n}\to \G$, $j_{n} : F_{n}\to \G$ such that $i_n\restriction E_{n-1}$ is $2^{-k}$-close to $i_{n-1}$, $j_n\restriction F_{n-1}$ is $2^{-k}$-close to $j_{n-1}$ and $U\cmp i_n$ is $2^{-k}$-close to $j_n\cmp T_n$.

Let $\{T_n:E_n\to F_n \}_{n\in \N}$ be the chain of non-expansive operators between finite-dimensional normed spaces resulting from a
fixed play, when Adam was using his strategy. In particular, Adam has recorded sequences $\{T_n:E_n\to F_n\}_{n\in \N}$, $\{i_n:E_n\to \G\}_{n\in \N}$,  $\{j_n:F_n\to \G\}_{n\in \N}$ of linear isometric embeddings such that $i_{2n+1}\restriction E_{2n-1}$ is $2^{-n}$-close to $i_{2n-1}$ and $j_{2n+1}\restriction F_{2n-1}$ is $2^{-n}$-close to $j_{2n-1}$ for each $n \in \N$. 

Let $\map{T_\infty}{E_\infty}{F_\infty}$ denote the completion of those unions, namely, $E_\infty$ is the completion of $\sett{E_n}{\ntr}$, $F_\infty$ is the completion of $\sett{F_n}{\ntr}$ and $T_\infty \rest E_n = T_n$ for every $\ntr$.
The assumptions that $i_{2n+1}[E_{2n+1}]$ contains all the vectors $a_0, \dots, a_n$ and $j_{2n+1}[F_{2n+1}]$ contains all the vectors $b_0, \dots, b_n$ ensures that both $i_\infty[E_\infty]$, $j_\infty[F_\infty]$ are dense in $\G$, where $\map{i_\infty}{E_\infty}{\G}$, $\map{j_\infty}{F_\infty}{\G}$ are pointwise limits of $\sett{i_n}{\Ntr}$ and $\sett{j_n}{\Ntr}$, respectively.
More precisely, $i_\infty \rest E_k$ is the pointwise limit of $\sett{i_n \rest E_k}{n \goe k}$ and $j_\infty \rest F_k$ is the pointwise limit of $\sett{j_n \rest F_k}{n \goe k}$ for every $k \in \Ntr$.
In particular, both $i_\infty$ and $j_\infty$ are surjective linear isometries.

Finally, $U \cmp i_\infty = j_\infty \cmp T_\infty$, because $U \cmp i_{2k}$ is $2^{-k}$-close to $j_{2k} \cmp T_{2k}$ for every $k \in \N$.
This completes the proof.

\end{pf}

\begin{pyt}\label{QUEbsdhkbds}
	Is $\uop$ generic in the space of all non-expansive operators on the \Gurarii\ space? Being ``generic" means of course that the set
	$$\setof{i \cmp \uop \cmp j}{i,j \text{ bijective linear isometries of }\G}$$
	is residual in the space of all non-expansive operators on $\G$. Here, it is natural to consider the pointwise convergence (i.e., strong operator) topology.
\end{pyt}

One could also consider a ``parametrized'' variant of the game above, where the two players build a chain of non-expansive operators from finite-dimensional normed spaces into a fixed Banach space $\bS$. If $\bS$ is separable then similar arguments as in the proof of Theorem~\ref{THMsztyryJedna} show that the second player has a strategy leading to $\upro{\bS}$. Thus, a variant of Question~\ref{QUEbsdhkbds} makes sense: Is it true that isometric copies of $\upro{\bS}$ form a residual set in a suitable space of operators?

After concluding that $\uop = \upro{\G}$, it seems that the ``parametrized'' construction of universal projections is better in the sense that it ``captures'' both the \Gurarii\ space $\G$ (when the range is the trivial space $\sn0$) and the universal operator $\uop$ (when the range equals $\G$), but also other examples, including projections from the \Gurarii\ space onto any separable Lindenstrauss space (see~\cite{Wojt} and~\cite{CsGwK}).

\separator

\paragraph{Acknowledgments.}
The authors would like to thank Przemys{\l}aw Wojtaszczyk for pointing out Example~\ref{ExUnivOptLeft}.

\end{document}